\documentclass[11pt,a4paper,twoside,fleqn]{article}
\usepackage[top=0.77in,left=1.58in,right=1.77in,bottom=3.52in]{geometry}
\usepackage{amsmath,amsfonts,amssymb}
\usepackage[none]{hyphenat}
\usepackage{graphicx}
\usepackage{mathptmx}
\usepackage[margin=10pt,font=small,labelfont=bf,labelsep=quad]{caption}
\usepackage{avant}
\usepackage{fancyhdr}
\usepackage{enumitem}
\usepackage[hang,flushmargin]{footmisc}
\usepackage[explicit]{titlesec}
\usepackage [running,mathlines] {lineno}

\usepackage{ITBjournal}

\journalID{
ITB J. ................... Vol. XX ..., No. X, 20XX, XX-XX}

\PublicationDate
{DD-MM-YYYY
}
{DD-MM-YYYY 
}
{DD-MM-YYYY 
}

\newcommand{\titlehead}{Inclusion Properties of Orlicz and Weak Orlicz Spaces}
\newcommand{\authorhead}{Al Azhary Masta, et.al 
}
\makeheader

\begin{document}
\newtheorem{theorem}{Theorem}[section]
\newtheorem{lemma}[theorem]{Lemma}
\newtheorem{corollary}[theorem]{Corollary}
\newtheorem{proposition}[theorem]{Proposition}
\newtheorem{conjecture}[theorem]{Conjecture}
\newtheorem{problem}[theorem]{Problem}
\newtheorem{claim}[theorem]{Claim}
\newtheorem{assumption}[theorem]{Assumption}
\newtheorem{remark}[theorem]{Remark}
\newtheorem{definition}[theorem]{Definition}
\newtheorem{example}[theorem]{Example}
\newtheorem{notation}{Notasi}
\renewcommand{\thenotation}{}


\begin{center}
\thispagestyle{fancy}
\title{Inclusion Properties of Orlicz and Weak Orlicz Spaces}
\author{Al Azhary Masta${}^{1,a}$, Hendra Gunawan${}^2$, Wono Setya Budhi${}^3$}
${}^{1,2,3}$Analysis and Geometry Group, Faculty of Mathematics
and Natural Sciences, Bandung Institute of Technology, Jl. Ganesha 10, Bandung 40132\\
${}^{a}$ \emph{Permanent Address}: Mathematics Education Department, Universitas Pendidikan Indonesia Jl. Dr. Setiabudi 229, Bandung 40154\\
Email: ${}^{1}$alazhari.masta@upi.edu, ${}^{2}$hgunawan@math.itb.ac.id, ${}^{3}$wono@math.itb.ac.id
\end{center}

\begin{abstract}
In this paper we discuss the structure of Orlicz spaces and weak
Orlicz spaces on $\mathbb{R}^n$. We obtain some necessary and sufficient
conditions for the inclusion property of these spaces. One of
the keys is to compute the norm of the characteristic functions
of the balls in $\mathbb{R}^n$.
\end{abstract}

\Keywords{Inclusion property, Orlicz spaces, weak Orlicz spaces.}

\section{Introduction}
Orlicz spaces were introduced by Z.W. Birnbaum and W. Orlicz in 1931
\cite{Orlicz}. Let $\Phi:[0,\infty)\to [0,\infty)$ be a Young function,
that is, $\Phi$ is
convex, left-continuous, $\Phi(0) = 0$, and $\lim \limits_{t\to\infty}
\Phi(t) = \infty$. Given a measure space $(X,dx)$, we define the Orlicz
space $L_\Phi(X)$ to be
the set of measurable functions $f : X \rightarrow \mathbb{R}$ such that
$$ \int_{X} \Phi ( a|f(x)|) dx < \infty$$
for some $a > 0$.  The space $L_\Phi(X)$ is a Banach space when equipped with
the norm
$$
\| f \|_{L_\Phi(X)} := \inf \left\{  {b>0:
\int_{X}\Phi \left(\frac{|f(x|)}{b} \right) dx \leq1}\right\}
$$
(see \cite{Christian, Luxemburg}). Note that, if $\Phi(t) := t^p$ for some $p\geq 1$ and $X:=\mathbb{R}^n$, then
$L_\Phi(X) = L_p(\mathbb{R}^n)$, the Lebesgue space of $p$-th integrable
functions on $\mathbb{R}^n$ \cite{Al}. Thus, Orlicz spaces can be viewed as a
generalization of Lebesgue spaces.

Many authors have been culminating important observations about Orlicz
spaces (see \cite{Kufner, Christian, Luxemburg, Al,Alen, Rao, Welland}, etc.) Here we are interested in the inclusion property of these
spaces. In \cite{Welland}, R. Welland proved the following
inclusion property: Let $X$ be of finite measure, and $\Phi,
\Psi$ be two Young functions. If there is $C>0$ such that $\Phi(t)\leq
\Psi(Ct)$ for every $t>0$, then $L_\Phi(X)\subseteq L_\Psi(X)$.
Accordingly, if $X$ is of finite measure, $\Phi, \Psi$ are two Young
functions, and there is $C>0$ such that $\Psi (\frac{t}{C})\leq \Phi(t)
\leq \Psi(Ct)$ for every $t>0$, then we have $L_\Phi(X)= L_\Psi(X)$.
A refinement of this result may be found in \cite{Kufner}, which states
that $L_\Phi(X)\subseteq L_\Psi(X)$ if only if there are $C>0$ and
$T > 0$ such that $\Phi(t)\leq\Psi(Ct)$ for every $ t\geq T $. Related results can be found in \cite{Abdullah,Mursaleen,Alen}. Motivated by these results, the purpose of this study is to get the inclusion property of Orlicz spaces $L_{\Phi}(\mathbb{R}^n)$ and
extend the results to weak Orlicz spaces $wL_{\Phi}(\mathbb{R}^n)$ (see \cite{Bekjan, Ning}).
(Here $X:=\mathbb{R}^n$ has an infinite measure.)

The rest of this paper is organized as follows. The main results are
presented in Sections 2 and 3. In Section 2, we state the inclusion
property of Orlicz spaces $L_{\Phi}(\mathbb{R}^n)$ as Theorem
\ref{theorem:2.5}, which contains a necessary and sufficient condition
for the inclusion property to hold. An analogous result for the weak
Orlicz spaces $wL_{\Phi}(\mathbb{R}^n)$ is stated as Theorem 3.3.

To prove the results, we pay attention to the characteristic functions
of balls in $\mathbb{R}^n$ and use the inverse function of $\Phi$, namely
$\Phi^{-1}(s) :=\inf \{ r \geq 0 : \Phi (r) > s \}$. The reader will
find the following lemma useful.

\medskip

\begin{lemma}\label{lemma:1.1}
Suppose that $\Phi$ is a Young function and $ \Phi^{-1}(s)=\inf \{r \geq 0 :
\Phi (r) > s \}$. We have

{\parindent=0cm
{\rm (1)} $\Phi^{-1}(0) = 0$.

{\rm (2)} $ \Phi^{-1}(s_1) \leq \Phi^{-1}(s_2)$ for  $s_1 \leq s_2$.

{\rm (3)} $\Phi (\Phi^{-1}(s)) \leq s \leq \Phi^{-1}(\Phi(s))$ for $0 \leq s <
\infty$.

{\rm (4)} Let $C>0$. Then $\Phi_1(t) \leq \Phi_2(Ct)$ if only if $C\Phi^{-1}_1(t)
\geq \Phi^{-1}_2(t)$, for every $ t \geq 0$.

{\rm (5)} Let $C>0$. Then $\Phi_1(t) \leq C \Phi_2(t)$ if only if  $\Phi^{-1}_1(Ct)
\geq \Phi^{-1}_2(t)$, for every $ t \geq 0$.
\par}
Proof.
\end{lemma}
The proof of parts (1)-(3) can be found in \cite{Nakai1}. Now we will prove (4) and (5).

{\parindent=0cm
{\rm (4)} Let $C>0$. We will prove $\Phi_1(t) \leq \Phi_2(Ct)$ if only if $C\Phi^{-1}_1(t)
\geq \Phi^{-1}_2(t)$, for every $ t \geq 0$.\\
Take an arbitrary $C > 0$ such that $\Phi_1(t) \leq \Phi_2(Ct)$ for every $ t \geq 0$. Let $ \Phi^{-1}_1(t) = \inf A_1$ for $ A_1 = \inf \{r \geq 0 : \Phi_{1} (r) > t \}$ and $B_1 =  \{r \geq 0 : \Phi_{2} (Cr) > t \}$. Observe that

\begin{align*}
\inf B_1 &= \inf \{r \geq 0 : \Phi_{2} (Cr) > t \}\\
      & = \inf \{\frac{x}{C} \geq 0 : \Phi_{2} (x) > t \}\\
      & = \frac{1}{C} \inf \{x \geq 0 : \Phi_{2} (x) > t \}\\
      & = \frac{1}{C} \Phi^{-1}_2(t)
\end{align*}
for $ x = Cr$. Take an arbitrary $ r \in A_1$, we have $ \Phi_2(Cr) \geq \Phi_1(r) > t $. Hence it follows that $ r \in B_1$, and so we conclude that $ A_1 \subseteq B_1$. Accordingly, we obtain $\frac{1}{C} \Phi^{-1}_2(t) = \inf B_1 \leq \inf A_1 = \Phi^{-1}_1(t)$. \\
Now, suppose that $C\Phi^{-1}_1(t) \geq \Phi^{-1}_2(t)$, for $ C > 0$ and every $ t \geq 0$. Observe that, by Lemma \ref{lemma:1.1} (3) we have
\begin{align*}
\Phi_{1}(\frac{t}{C}) & \leq \Phi_{1}(\frac{\Phi^{-1}_{2}(\Phi_2(t))}{C})\\
& \leq \Phi_{1}(C \frac{\Phi^{-1}_{1}(\Phi_2(t))}{C})\\
& = \Phi_{1}(\Phi^{-1}_{1}(\Phi_2(t)))\\
& \leq \Phi_{2}(t).
\end{align*}

As a result, we have $\Phi_{1}(\frac{t}{C})\leq \Phi_{2}(t)$ or $ \Phi_1(t) \leq \Phi_2(Ct)$.

{\rm (5)} Let $C>0$. We will prove $\Phi_1(t) \leq C \Phi_2(t)$ if only if  $\Phi^{-1}_1(Ct)
\geq \Phi^{-1}_2(t)$, for every $ t \geq 0$.
\par}
Take an arbitrary $C > 0$ such that $\Phi_1(t) \leq C \Phi_2(t)$ for every $ t \geq 0$. Let $ \Phi^{-1}_1(Ct) = \inf A_2$ for $ A_2 = \inf \{r \geq 0 : \Phi_{1} (r) > Ct \}$ and $\Phi^{-1}_2(t) = \inf B_2$ for $B_2 =  \{r \geq 0 : \Phi_{2} (r) > t \}$. Observe that, for $r \in A_2$ we have $C \Phi_2(t) \geq \Phi_{1} (r) > Ct$. Hence it follows that $ r \in B_2$, and so we conclude that $ A_2 \subseteq B_2$. Accordingly, we obtain $\Phi^{-1}_2(t) = \inf B_2 \leq\inf A_2 = \Phi^{-1}_1(Ct)$.\\
Now, suppose that $\Phi^{-1}_1(Ct) \geq \Phi^{-1}_2(t)$, for $ C > 0$ and every $ t \geq 0$. Observe that, by by Lemma \ref{lemma:1.1} (3) we have
\begin{align*}
\Phi_{1}(t) & \leq \Phi_{1}(\Phi^{-1}_{2}(\Phi_2(t)))\\
& \leq \Phi_{1}(\Phi^{-1}_{1}(C\Phi_2(t))) \\
& \leq C\Phi_{2}(t).
\end{align*}
As a result, we have $\Phi_{1}(t)\leq C\Phi_{2}(t)$.\\
More lemmas (and their proofs) will be presented in the next sections.

\section{Inclusion property of Orlicz spaces}
Let us first recall several Lemmas in the following.

\begin{lemma}\label{lemma:2.1}
Let $\Phi$ be a Young function, then $\Phi (\alpha t) \leq \alpha \Phi (t)$ for $t>0$ and $0\leq \alpha \leq 1$.
\end{lemma}

\begin{lemma}\label{lemma:2.2} \cite{Christian}
 Let $\Phi$ be a Young function and $f \in L_{\Phi(\mathbb{R}^n)}$. If $0 <
 \| f \|_{L_\Phi(\mathbb{R}^n)} <\infty$, then $\int_{\mathbb{R}^n}\Phi\left(
 \frac{|f(x)|}{\| f \|_{L_\Phi(\mathbb{R}^n)}} \right)  dx \leq 1 $. Furthermore,
 $\| f \|_{L_\Phi(\mathbb{R}^n)} \leq 1 $ if only if $\int_{\mathbb{R}^n}\Phi(|f(x)|) dx
 \leq 1.$
\end{lemma}

\begin{corollary}\label{corollary:2.3}
Let $\Phi, \Psi$ be Young functions. If there exists $C > 0$ such that
$\Phi(t) \leq \Psi(Ct)$ for $ t > 0$, then
$L_{\Psi}(\mathbb{R}^n) \subseteq L_{\Phi}(\mathbb{R}^n)$ with $\| f \|_{L_\Phi(\mathbb{R}^n)} \leq
C \| f \|_{L_\Psi(\mathbb{R}^n)}$
for every $f\in L_{\Psi}(\mathbb{R}^n)$.
\end{corollary}

\noindent{\it Proof}. Suppose that $f\in L_\Psi(\mathbb{R}^n)$.
Observe that
\[
\int_{\mathbb{R}^n}\Phi \left( \frac{|f(x)|}{C \| f \|_{L_\Psi(\mathbb{R}^n)}} \right) dx
\leq \int_{\mathbb{R}^n}\Psi \left(\frac{C|f(x)|}{ C\| f \|_{L_\Psi(\mathbb{R}^n)}} \right)dx=
\int_{\mathbb{R}^n} \Psi \left( \frac{|f(x)|}{ \| f \|_{L_\Psi(\mathbb{R}^n)}} \right) dx\le 1.
\]
By the definition of $\| f \|_{L_\Phi(\mathbb{R}^n)}$, we have $\| f \|_{L_\Phi(\mathbb{R}^n)}
\leq C \| f \|_{L_\Psi(\mathbb{R}^n)}.$
This proves that $L_{\Psi}(\mathbb{R}^n) \subseteq L_{\Phi}(\mathbb{R}^n)$, as desired.

\noindent{\tt Remark}. From Corollary \ref{corollary:2.3}, we note that
if $\Phi \le \Psi$, then $L_\Psi(\mathbb{R}^n) \subseteq L_\Phi(\mathbb{R}^n)$ with
$\| f \|_{L_\Phi(\mathbb{R}^n)}\leq \| f \|_{L_\Psi(\mathbb{R}^n)}$ for every $f\in
L_\Psi(\mathbb{R}^n)$. As we shall see below, the converse of this statement also
holds. We need the following lemma.

\begin{lemma}\label{lemma:2.4}\cite{Rao} Let $\Phi$ be a Young function,
$a\in\mathbb{R}^n$, and $r>0$. Then
$ \| \chi_{B(a,r)} \|_{L_\Phi(\mathbb{R}^n)} = \frac{1}{\Phi^{-1}(\frac{1}
{|B(a,r)|})}$, where $|B(a,r)|$ denotes the volume of $B(a,r)$.
\end{lemma}

\begin{theorem}\label{theorem:2.5}
Let $\Phi, \Psi$ be Young functions. Then the following statements are equivalent:

{\parindent=0cm
{\rm (1)} $\Phi(t) \leq \Psi(Ct)$ for every $t > 0.$

{\rm (2)} $L_{\Psi}(\mathbb{R}^n) \subseteq L_{\Phi}(\mathbb{R}^n)$.

{\rm (3)} For every $ f \in L_{\Psi}(\mathbb{R}^n)$, we have  $\| f \|_{L_\Phi(\mathbb{R}^n)} \leq
  C \| f \|_{L_\Psi(\mathbb{R}^n)}.$
\par}
\end{theorem}

\noindent{\it Proof}.
We have seen that (1) implies (2). Next, since $(L_{\Psi}(\mathbb{R}^n),L_{\Phi}(\mathbb{R}^n))$ is a Banach pair, it follows from  [5, Lemma 3.3] that (2) and (3) are equivalent. It thus remains to show that (3) implies (1). Now assume that (3) holds. By Lemma
\ref{lemma:2.4}, we have
$$
\frac{1}{\Phi^{-1}(\frac{1}{|B(a,r)|})} = \| \chi_{B(a,r)} \|_{L_\Phi(\mathbb{R}^n)}
\leq C\| \chi_{B(a,r)} \|_{L_\Psi(\mathbb{R}^n)} =  \frac{C}{\Psi^{-1}(\frac{1}{|B(a,r)|})},
$$
or $ C\Phi^{-1}(\frac{1}{|B(a,r)|}) \geq \Psi^{-1}(\frac{1}{|B(a,r)|})$,
for every $a\in\mathbb{R}^n,\ r>0$.
By Lemma \ref{lemma:1.1}(4), we obtain $\Phi(\frac{1}{|B(a,r)|}) \leq
\Psi(\frac{C}{|B(a,r)|})$. Since $ r > 0$ is arbitrary, we conclude that
$\Phi(t) \leq \Psi(Ct)$ for every $t > 0.$

\subsection{A special case} 
One may ask whether from inclusion relations between Orlicz spaces, we may deduce some known fact of those in Lebesgue spaces. 
The answer is affirmative; we need the following lemma for the purpose.

\begin{lemma}\label{lemma:2.6}
Let $\Phi_1, \Phi_2$, and $\Phi_3$ be Young functions such that
$ \Phi^{-1}_1(t)\Phi^{-1}_2(t) \leq \Phi^{-1}_3(t)$ for every $ t \geq 0$.
If $f \in L_{\Phi_1}(\mathbb{R}^n)$ and $g \in L_{\Phi_2}(\mathbb{R}^n)$, then
$fg \in L_{\Phi_3}(\mathbb{R}^n)$ with $\| fg \|_{L_{\Phi_3}(\mathbb{R}^n)}
\leq 2 \| f \|_{L_{\Phi_1}(\mathbb{R}^n)} \| g \|_{L_{\Phi_2}(\mathbb{R}^n)}$.
\end{lemma}

\noindent{\it Proof}. Let $s,t \ge 0$.
Without loss of generality, suppose that $\Phi_1(s) \leq \Phi_2(t)$. By Lemma \ref{lemma:1.1}(3), we obtain
$$
st \leq \Phi^{-1}_1(\Phi_1(s))\Phi^{-1}_2(\Phi_2(t))\leq \Phi^{-1}_1(\Phi_2(t))
\Phi^{-1}_2(\Phi_2(t))\leq \Phi^{-1}_3(\Phi_2(t)).
$$
Hence $\Phi_3(st)\leq \Phi_3(\Phi^{-1}_3(\Phi_2(t))) \leq \Phi_2(t) \leq
\Phi_2(t) + \Phi_1(s)$. From Lemma \ref{lemma:2.1}, we have

$$\int_{\mathbb{R}^n}\Phi_3 \left( \frac{|f(x)g(x)|}{2 \| f \|_{L_{\Phi_1}(\mathbb{R}^n)}
\| g \|_{L_{\Phi_2}(\mathbb{R}^n)}} \right) dx \leq \frac{1}{2} \int_{\mathbb{R}^n}\Phi_3 \left(
\frac{|f(x)g(x)|}{\| f \|_{L_{\Phi_1}(\mathbb{R}^n)} \| g \|_{L_{\Phi_2}(\mathbb{R}^n)}}\right) dx.$$

On the other hand, by Lemma \ref{lemma:2.2} we obtain\\
$\int_{\mathbb{R}^n}\Phi_3 \left(
\frac{|f(x)g(x)|}{\| f \|_{L_{\Phi_1}(\mathbb{R}^n)} \| g \|_{L_{\Phi_2}(\mathbb{R}^n)}}\right) dx \leq \int_{\mathbb{R}^n}\Phi_1 \left( \frac{|f(x)|}{\| f \|_{L_{\Phi_1}(\mathbb{R}^n)}}
\right) dx + \int_{\mathbb{R}^n}\Phi_2 \left( \frac{|g(x)|}{\| g \|_{L_{\Phi_2}(\mathbb{R}^n)}}
\right) dx \leq 2,$

whenever $f \in L_{\Phi_1}(\mathbb{R}^n)$ and $ g \in L_{\Phi_2}(\mathbb{R}^n)$. By using the
definition of $\| fg \|_{L_{\Phi_3}(\mathbb{R}^n)}$, we have $\| fg \|_{L_{\Phi_3}(\mathbb{R}^n)} \leq
2 \| f \|_{L_{\Phi_1}(\mathbb{R}^n)}\|  g \|_{L_{\Phi_2}(\mathbb{R}^n)}$, as desired.

\begin{corollary}\label{corollary:2.7}
Let $X:=B(a,r_0) \subset\mathbb{R}^n$ for some $a\in\mathbb{R}^n$ and $r_{0} > 0$.
If $\Phi_1, \Phi_2$ are
two Young functions and there is a Young function $\Phi$ such that
$$
\Phi^{-1}_1(t)\Phi^{-1}(t) \leq \Phi^{-1}_2(t)
$$
for every $ t \geq 0$, then $L_{\Phi_1}(X)\subseteq L_{\Phi_2}(X)$ with
$$
\| f \|_{L_{\Phi_2}(X)} \leq \frac{2}{\Phi^{-1}(\frac{1}
{|B(a,r_{0})|})} \| f \|_{L_{\Phi_1}(X)}
$$
for $ f \in L_{\Phi_1}(X)$.
\end{corollary}

\noindent{\it Proof}. Let $f \in L_{\Phi_1}(X)$.
By Lemma \ref{lemma:2.4} and choosing $g:=\chi_{B(a,r_0)}$,
we obtain
$$
\| f \chi_{B(a,r_0)} \|_{L_{\Phi_2}(X)} \leq 2 \| \chi_{B(a,r_0)}
\|_{L_\Phi(X)} \| f \|_{L_{\Phi_1}(X)} = \frac{2}{\Phi^{-1}(\frac{1}{|B(a,r_{0})|})}
\| f \|_{L_{\Phi_1}(X)}.
$$
This shows that $L_{\Phi_1}(X) \subseteq L_{\Phi_2}(X)$.

\begin{corollary}
Let $X:=B(a,r_0)$ for some $a\in\mathbb{R}^n$ and $r_0>0$. If $ 1 \leq p_{2} <
p_{1} < \infty$, then $L_{p_{1}}(X) \subseteq L_{p_{2}}(X)$.
\end{corollary}

\noindent{\it Proof}.
Let $\Phi_1(t):= t^{p_1}, \Phi_2(t):= t^{p_2}$, and $\Phi(t):=
t^{\frac{p_1 p_2}{p_1- p_2}}$ ($t\ge 0$). Since $ 1 \leq p_{2} < p_{1} < \infty$,
we have $\frac{p_1 p_2}{p_1- p_2} > 1$. Thus, $\Phi_1,\
\Phi_2$, and $\Phi$ are three Young functions. Observe that, using the definition
of $\Phi^{-1}$ and Lemma \ref{lemma:1.1}, we have
$$
\Phi^{-1}_{1}(t) =t^{\frac{1}{p_1}},\ \Phi^{-1}_{2}(t) =t^{\frac{1}{p_2}},\
{\rm and}\ \Phi^{-1}(t) =t^{\frac{p_1- p_2}{p_1p_2}}.
$$
Moreover, $\Phi^{-1}_{1}(t)\Phi^{-1}(t) = t^{\frac{1}{p_1}}t^{\frac{p_1-
p_2}{p_1p_2}} = t^{\frac{1}{p_2}} = \Phi^{-1}_{2}(t)$, and so it follows
from Corollary \ref{corollary:2.7} that $\| f \|_{L_{p_2}(X)} \leq \frac{2}{\Phi^{-1}(
\frac{1}{|B(a,r_0)|})}  \| f \|_{L_{p_1}(X)}$, and therefore $L_{p_{1}}(X)
\subseteq L_{p_{2} }(X)$.

\noindent{\tt Remark}. Of course we can prove the inclusion of property of
Lebesgue spaces on a finite measure space directly via H\"older's inequality.
What we showed here is that we can obtain the result through the lens of
Orlicz spaces.

\section{Inclusion property of Weak Orlicz spaces}

First, we recall the definition of weak Orlicz spaces \cite{Nakai2}.
Let $\Phi$ be a Young function. We define the weak Orlicz spaces $wL_{\Phi}(\mathbb{R}^n)$
to be the set of measurable functions $f : \mathbb{R}^n \rightarrow \mathbb{R} $ such
that $\| f \|_{wL_\Phi(\mathbb{R}^n)} < \infty$, where
$$
\| f \|_{wL_\Phi(\mathbb{R}^n)} := \inf \left\{  {b>0:
\mathop {\sup }\limits_{t > 0} \Phi(t) \Bigl| \{ x \in \mathbb{R}^n : \frac{|f(x)|}{b} > t \}
\Bigr| \leq1}\right\}.
$$

\noindent{\tt Remark}. Note that $\| \cdot \|_{wL_\Phi(\mathbb{R}^n)}$ defines a quasi-norm in $wL_{\Phi}(\mathbb{R}^n)$, and that \; \; \; \; \;    $(wL_{\Phi}(\mathbb{R}^n),\| \cdot \|_{wL_\Phi(\mathbb{R}^n)})$ forms a quasi-Banach space (see \cite{Bekjan, Yong}).

 The relation between weak Orlicz spaces and (strong) Orlicz spaces is clear, as
presented in the following theorem.

\begin{theorem}\label{theorem:3.1}\cite{Yong,Xueying}
Let $\Phi$ be a Young function. Then $L_{\Phi}(\mathbb{R}^n) \subset wL_{\Phi}(\mathbb{R}^n)$
with $\| f \|_{wL_\Phi(\mathbb{R}^n)} \leq \| f \|_{L_\Phi(\mathbb{R}^n)}$ for every $f\in L_\Phi(\mathbb{R}^n)$.
\end{theorem}

\noindent{\it Proof}. The proof of this theorem can be found in \cite{Yong,Xueying}. We rewrite the proof here for convenience. 

Given $f\in L_\Phi(\mathbb{R}^n)$,
let $A_{\Phi,w} := \Bigl\{  {b>0: \mathop {\sup }\limits_{t > 0} \Phi(t)
\bigl| \{ x \in \mathbb{R}^n : \frac{|f(x)|}{b} > t \} \bigr| \leq 1}\Bigr\} $ and $ B_{\Phi,w} :=
\left\{  {b>0: \int_{\mathbb{R}^n}\Phi \left(\frac{|f(x|)}{b} \right) dx \leq1}\right\}$.
Then $\| f \|_{wL_\Phi(\mathbb{R}^n)} = \inf A_{\Phi,w}$ and $\| f \|_{L_\Phi(\mathbb{R}^n)} = \inf B_{\Phi,w}$.
Observe that, for arbitrary $ b \in B_{\Phi,w}$ and $t > 0$, we have
\begin{align*}
\Phi(t) \Bigl| \{ x \in \mathbb{R}^n : \frac{|f(x)|}{b} > t \} \Bigr|
 & \leq
\int_{\{ x \in \mathbb{R}^n : \frac{|f(x)|}{b} > t \} }\Phi \left(\frac{|f(x)|}{b} \right) dx \\
 & \leq
\int_{\mathbb{R}^n}\Phi \left(\frac{|f(x)|}{b} \right) dx  \leq 1.
\end{align*}
Since $ t > 0$ is arbitrary, we have $\mathop {\sup }\limits_{t > 0} \Phi(t)
\mid \{ x \in \mathbb{R}^n : \frac{|f(x)|}{b} > t \} \mid \leq1$, and
$ B_{\Phi,w} \subseteq A_{\Phi,w}$. Hence, $f\in wL_{\Phi}$ with
$\| f \|_{wL_\Phi(\mathbb{R}^n)} \leq \| f \|_{L_\Phi(\mathbb{R}^n)}$.

\noindent{\tt Remark}. As the strong and weak Orlicz spaces contain the
strong and weak Lebesgue spaces respectively, the inclusion in the above
theorem is proper. See \cite{Gunawan} for a counterexample.

In addition to Lemma \ref{lemma:2.4}, we have the following lemma for the characteristic
functions of balls in weak Orlicz spaces.

\begin{lemma}\label{lemma:3.2}\cite{Ning}
Let $\Phi$ be a Young function, $a\in\mathbb{R}^n$, and $r >0$ be arbitrary. Then we have
$ \| \chi_{B(a,r)} \|_{wL_\Phi(\mathbb{R}^n)} = \frac{1}{\Phi^{-1}(\frac{1}{|B(a,r)|})}$.
\end{lemma}

Now we come to the inclusion property of weak Orlicz spaces.
\begin{theorem}\label{theorem:3.3}
Let $\Phi, \Psi$ be Young functions. Then the following statements are equivalent:

{\parindent=0cm
{\rm (1)} $\Phi(t) \leq \Psi(Ct)$ for every $t > 0$.

{\rm (2)} $wL_{\Psi}(\mathbb{R}^n) \subseteq wL_{\Phi}(\mathbb{R}^n)$.

{\rm (3)} For every $ f \in wL_{\Psi}(\mathbb{R}^n)$, we have $\| f \|_{wL_\Phi(\mathbb{R}^n)} \leq C\| f \|_{wL_\Psi(\mathbb{R}^n)}.$

\par}
\end{theorem}

\noindent{\it Proof}.
Assume that (1) holds, and let $ f \in wL_{\Psi}(\mathbb{R}^n)$.
Put
$$A_{\Phi,w} = \Bigl\{  {b>0: \mathop {\sup }\limits_{t > 0} \Phi(t) \bigl|
\{ x \in \mathbb{R}^n : \frac{|f(x)|}{b} > t \} \bigr| \leq1}\Bigr\}
$$
and
\begin{align*}
A_{\Psi,w} & = \Bigl\{  {b>0: \mathop {\sup }\limits_{t > 0} \Psi(Ct) \bigl| \{ x \in \mathbb{R}^n : \frac{|f(x)|}{b} > t \} \bigr| \leq1}\Bigr\} \\
& = \Bigl\{  {b>0: \mathop {\sup }\limits_{s > 0} \Psi(s) \bigl|
\{ x \in \mathbb{R}^n : \frac{C|f(x)|}{b} > s \} \bigr| \leq1}\Bigr\},
\end{align*}

for $s=Ct$. Then $\| f \|_{wL_\Phi(\mathbb{R}^n)} = \inf A_{\Phi,w}$  and $C\| f \|_{wL_\Psi(\mathbb{R}^n)} = \| Cf \|_{wL_\Psi(\mathbb{R}^n)} = \inf A_{\Psi,w}$.
Observe that, for arbitrary $ b \in A_{\Psi,w} $ and $ t > 0$, we have
$$
\Phi(t) \bigl| \{ x \in \mathbb{R}^n : \frac{|f(x)|}{b} > t \} \bigr| \leq  \Psi(Ct) \bigl|
\{ x \in \mathbb{R}^n : \frac{|f(x)|}{b} > t \} \bigr| \leq 1.
$$
Thus, $\mathop {\sup }\limits_{t > 0} \Phi(t) \bigl| \{ x \in \mathbb{R}^n :
\frac{|f(x)|}{b} > t \} \bigr| \leq 1$. Hence it follows that $ b \in A_{\Phi,w}$,
and so we conclude that $A_{\Psi,w} \subseteq A_{\Phi,w}$. Accordingly, we obtain
$$\| f \|_{wL_\Phi(\mathbb{R}^n)} = \inf A_{\Phi,w} \leq \inf A_{\Psi,w} = C \| f \|_{wL_\Psi(\mathbb{R}^n)},$$
which also proves that $wL_{\Psi}(\mathbb{R}^n) \subset wL_{\Phi}(\mathbb{R}^n)$.

As mentioned in [10, Appendix G], we are aware that [5, Lemma 3.3] still holds for quasi-Banach spaces, and so (2) and (3) are equivalent.

Assume now that (3) holds. By Lemma \ref{lemma:3.2}, we have
$$
 \frac{1}{\Phi^{-1}(\frac{1}{|B(a,r_{0})|})} = \| \chi_{B(a,r_{0})} \|_{wL_\Phi(\mathbb{R}^n)}
 \leq C \| \chi_{B(a,r_{0})} \|_{wL_\Psi(\mathbb{R}^n)} =  \frac{C}{\Psi^{-1}(\frac{1}{|B(a,r_{0})|})},
$$
or $ C \Phi^{-1}(\frac{1}{|B(a,r_{0})|}) \geq \Psi^{-1}(\frac{1}{|B(a,r_{0})|})$,
for arbitrary $a\in\mathbb{R}^n$ and $r_0 > 0$. By Lemma \ref{lemma:1.1}, we have
$$
\Phi\Bigl(\frac{1}{|B(a,r_{0})|}\Bigr) \leq  \Psi\Bigl(\frac{C}{|B(a,r_{0})|}\Bigr).
$$
Since $a\in\mathbb{R}^n$ and $r_{0} > 0$ are arbitrary, we conclude that $\Phi(t) \leq \Psi(Ct)$ for every $t > 0$.

\section{Concluding remarks}

We have shown the inclusion property of (strong) Orlicz spaces and of weak Orlicz
spaces. Both use the norm of the characteristic functions of the balls in $\mathbb{R}^n$.
As our final conclusion, we have the following corollary which states that the
inclusion property of (strong) Orlicz spaces are equivalent to that of weak
Orlicz spaces, and both can be observed just by comparing the associated Young
functions.\\
As a corollary of Theorem \ref{theorem:2.5} and Theorem \ref{theorem:3.3} we have that if $\Phi, \Psi$ are two Young functions, then the following statements are equivalent:

{\parindent=0cm
{\rm (1)} $\Phi(t) \leq \Psi(Ct)$ for every $t > 0$.

{\rm (2)} $L_{\Psi}(\mathbb{R}^n) \subseteq L_{\Phi}(\mathbb{R}^n).$

{\rm (3)} For every $ f \in L_{\Psi}(\mathbb{R}^n)$, we have $\| f \|_{L_\Phi(\mathbb{R}^n)}
\leq C\| f \|_{L_\Psi(\mathbb{R}^n)}$.

{\rm (4)} $wL_{\Psi}(\mathbb{R}^n) \subseteq wL_{\Phi}(\mathbb{R}^n)$.

{\rm (5)} For every $ f \in wL_{\Psi}(\mathbb{R}^n)$, we have $\| f \|_{wL_\Phi(\mathbb{R}^n)}
\leq \| f \|_{wL_\Psi(\mathbb{R}^n)}$.
\par}

\medskip

\section*{Acknowledgement}
The first and second authors are supported by ITB
Research and Innovation Program 2015.



\begin{thebibliography}{9}
\parskip=-5pt
\labelsep=7pt
\bibitem{Abdullah}
A. Alotaibi, M. Mursaleen, and S. K. Sharma, "Double sequence spaces over \emph{n}-normed spaces defined by a sequence of Orlicz functions", \emph{J. Inequal. Appl.} \textbf{2014}-1(2014),1--12.

\bibitem{Bekjan}

T.N. Bekjan, Z. Chen, P. Liu, and Y. Jiao,"Noncommutative weak Orlicz spaces and martingale inequalities", \emph{Studia Math.} \textbf{204}-3 (2011), 195--212.

\bibitem{Gunawan}
	H. Gunawan, D.I. Hakim, K.M. Limanta, A.A. Masta, "Inclusion properties of generalized Morrey spaces", to appear in Math. Nachr.

\bibitem{Yong}
Y. Jiao,"Embeddings between weak Orlicz martingale spaces", \emph{J. Math. Appl.} \textbf{378 } (2011), 220--229.  
               
\bibitem {Krein}
S.G Kre\v{i}n, Yu.\={I} Petun\={i}n, and E.M. Sem\"{e}nov, "\emph{Interpolation of Linear Operators}", Translation of Mathematical Monograph 54, American Mathematical Society, Providence, R.I., 1982.

\bibitem{Kufner}
    O. Kufner, O. John, and S. Fu\"{c}ik, \emph{Function Space}, Noordhoff International
    Publishing, Czechoslovakia, 1977.

\bibitem{Christian}
    C. L\'{e}onard, "Orlicz spaces", preprint. [http://cmap.polytechnique.fr/$\sim$leonard/papers/orlicz.pdf, accessed on August 17, 2015.]

\bibitem{Ning}
    N. Liu and Y. Ye, "Weak Orlicz space and its convergence theorems", \emph{Acta Math.
    Sci. Ser. B} {\bf 30}-5 (2010), 1492--1500.

\bibitem{Luxemburg}
    W.A.J. Luxemburg, \emph{Banach Function Spaces}, Thesis, Technische Hogeschool te Delft,
    1955.

\bibitem {Al}
    A.A. Masta, \emph{On Uniform Orlicz Spaces}, Thesis, Bandung Institute of Technology, 2013.

\bibitem{Mursaleen}
M. Mursaleen , Kuldip Raj, and Sunil K. Sharma, "Some spaces of difference sequences and lacunary statitical convergence in \emph{n}-normed space defined by sequence of Orlicz functions", \emph{Miskolc Mathematical Notes} \textbf{16}-1 (2015), 283--304.
	
\bibitem{Nakai1}
    E. Nakai, "On Orlicz-Morrey spaces", research report [http://repository.kulib.kyoto-u.ac.jp/dspace/bitstream/2433/58769/1/1520-
10.pdf, accessed on August 17, 2015.]

\bibitem{Nakai2}
    E. Nakai, "Orlicz-Morrey spaces and some integral operators", research report [http://repository.kulib.kyoto-u.ac.jp/dspace/bitstream/2433/26035/1/
1399-13.pdf, accessed on August 17, 2015.]

\bibitem{Orlicz}
    W. Orlicz, \emph{Linear Functional Analysis (Series in Real Analysis Volume 4)},
    World Scientific, Singapore, 1992.

\bibitem{Alen}
A. Osan\c{c}liol,"Inclusion between weighted Orlicz spaces", \emph{J. Inequal. Appl.} \textbf{2014}-390 (2014), 1--8.

\bibitem{Rao}
    M.M. Rao and Z.D. Ren, \emph{Theory of Orlicz spaces, volume 146 of Monographs and Textbooks in Pure and Applied Mathematics}, Marcel Dekker, Inc., New York, 1991.

\bibitem{Welland}
    R. Welland, "Inclusion relations among Orlicz spaces", \emph{Proc. Amer. Math. Soc}.
    {\bf 17}-1 (1966), 135--139.

\bibitem{Xueying}
Xueying Zhang and Chuanzhou Zhang, "Weak Orlicz spaces generated by concave functions", International Conference on Information Science and Technology (ICIST) 2011, 42--44.

\end{thebibliography}
\end{document}